\documentclass[11pt]{article}
\usepackage{latexsym}
\usepackage{amsfonts}
\usepackage{mathrsfs,amsthm,amsmath}
\usepackage{color}
\usepackage[round]{natbib}
\newtheorem{theorem}{Theorem}[section]

\newtheorem{condition}{Condition}

\newtheorem{proposition}[theorem]{Proposition}
\newtheorem{remark}{Remark}[section]

\begin{document}
\title{Asymptotic results for linear combinations of spacings generated by
i.i.d. exponential random variables}
\author{Camilla Cal\`i\thanks{Address: Dipartimento di Biologia,
Universit\`a di Napoli Federico II, Via Cintia, Complesso Monte S.
Angelo, 80126 Naples, Italy. e-mail: \texttt{camilla.cali@unina.it}}\and
Maria Longobardi\thanks{Address: Dipartimento di Biologia,
Universit\`a di Napoli Federico II, Via Cintia,
Complesso Monte S. Angelo, 80126 Naples, Italy. e-mail:
\texttt{maria.longobardi@unina.it}}\and
Claudio Macci\thanks{Address: Dipartimento di Matematica,
Universit\`a di Roma Tor Vergata, Via della Ricerca Scientifica,
I-00133 Roma, Italia. E-mail: \texttt{macci@mat.uniroma2.it}}\and
Barbara Pacchiarotti\thanks{Address: Dipartimento di Matematica,
Universit\`a di Roma Tor Vergata, Via della Ricerca Scientifica,
I-00133 Roma, Italia. E-mail: \texttt{pacchiar@mat.uniroma2.it}}}
\maketitle
\begin{abstract}
We prove large (and moderate) deviations for a class of linear combinations of
spacings generated by i.i.d. exponentially distributed random variables. We
allow a wide class of coefficients which can be expressed in terms of continuous
functions defined on $[0,1]$ which satisfy some suitable conditions. In this way
we generalize some recent results by \cite{GiulianoMacciPacchiarottiJSPI2015}
which concern the empirical cumulative entropies defined in
\cite{DicrescenzoLongobardiJSPI2009}.
\end{abstract}

\section{Introduction}
Empirical processes and their applications to statistics are widely studied (see e.g.
\cite{ShorackWellner} as a monograph on this topic). An important part of the results on this topic
concerns linear combinations of order statistics (called $L$-statistics) and, more in particular, linear
combinations of spacings (a spacing is a difference between two consecutive order statistics). Among the
references with results on large deviations for $L$-statistics here we recall \cite{Aleshkyavichene},
\cite{BentkusZikitis}, \cite{GroeneboomOosterhoofRuymgaart} and \cite{GroeneboomShorack}. In some cases
the large deviation results are formulated in terms of the concept of \emph{large deviation principle}
(see e.g. \cite{DemboZeitouni}) and, among the references with this kind of results, here we recall
\cite{Boistard} and \cite{DuffyMacciTorrisi}.

The aim of this paper is to generalize the results in \cite{GiulianoMacciPacchiarottiJSPI2015} concerning
a particular sequence of linear combinations of spacings $\{C_n:n\geq 1\}$ generated by a sequence of
independent and identically distributed (i.i.d. for short) random variables $\{X_n:n\geq 1\}$. We recall
that the random variables $\{C_n:n\geq 1\}$ are the empirical cumulative entropies defined in
\cite{DicrescenzoLongobardiJSPI2009} for a sequence of i.i.d. positive random variables $\{X_n:n\geq 1\}$
with a (common) absolutely continuous distribution function. Moreover the results in
\cite{GiulianoMacciPacchiarottiJSPI2015} concern the case of exponentially distributed random variables
$\{X_n:n\geq 1\}$ and, in such a case, the joint distribution of the spacings has some nice properties.
In this paper the random variables $\{X_n:n\geq 1\}$ are again exponentially distributed, and we allow a
wide class of sequences of linear combinations of spacings $\{C_n(w):n\geq 1\}$, where $w$ is a continuous
function on $[0,1]$ which satisfies some suitable conditions.

We conclude with the outline of the paper. Section \ref{sec:preliminaries} is devoted to some preliminaries;
in particular we also illustrate the connections with some references as
\cite{DicrescenzoLongobardiJSPI2009}, \cite{DicrescenzoLongobardiPROCEEDINGS} and \cite{GaoZhaoAS2011}.
In Section \ref{sec:results} we generalize the results in \cite{GiulianoMacciPacchiarottiJSPI2015}.
The connections between our moderate deviation result and the moderate deviation result for $L$-statistics in
\cite{GaoZhaoAS2011} is discussed in Section \ref{sec:connections-with-Gao-Zhao}. Finally, in Section
\ref{sec:empirical-estimators}, we discuss some possible choices of the function $w$ based on some empirical
entropies in the literature.

\section{Preliminaries}\label{sec:preliminaries}
We start with some preliminaries on large deviations. We also present
the sequence studied in this paper, and some connection with the literature.

\subsection{Preliminaries on large deviations}\label{sub:LD-preliminaries}
Here we briefly recall some basic preliminaries on large
deviations (see e.g. \cite{DemboZeitouni}, pages 4-5). Let
$\mathcal{X}$ be a topological space equipped with its completed
Borel $\sigma$-field. A sequence of $\mathcal{X}$-valued random
variables $\{Z_n:n\geq 1\}$ satisfies the \emph{large deviation
principle} (LDP for short) with speed function $v_n$ and rate
function $I$ if: $\lim_{n\to\infty}v_n=\infty$; the function
$I:\mathcal{X}\to[0,\infty]$ is lower semi-continuous; we have the
upper bound
$$\limsup_{n\to\infty}\frac{1}{v_n}\log P(Z_n\in C)\leq-\inf_{x\in C}I(x)\ \textrm{for all closed sets}\ C,$$
and the lower bound
$$\liminf_{n\to\infty}\frac{1}{v_n}\log P(Z_n\in O)\geq-\inf_{x\in O}I(x)\ \textrm{for all open sets}\ O.$$ 
A rate function $I$ is said to be \emph{good} if its level sets
$\{\{x\in\mathcal{X}:I(x)\leq\eta\}:\eta\geq 0\}$ are compact. In
the LDPs presented in this paper we always have
$\mathcal{X}=\mathbb{R}$. In some cases we apply the G\"{a}rtner
Ellis Theorem (see e.g. Theorem 2.3.6 in \cite{DemboZeitouni})
with the speed function $v_n$, and we obtain LDPs with good rate functions 
(see Propositions \ref{prop:LDP} and \ref{prop:MD}). Here we briefly recall
the statement of this theorem for real valued random variables: if there exists
$$\Lambda(\theta):=\lim_{n\to\infty}\frac{1}{v_n}\log\mathbb{E}[e^{v_n\theta Z_n}]\ \mbox{for all}\ \theta\in\mathbb{R},$$
the origin belongs to the interior of
$$\mathcal{D}(\Lambda):=\{\theta\in\mathbb{R}:\Lambda(\theta)<\infty\},$$
and the function $\Lambda$ is essentially smooth (see e.g.
Definition 2.3.5 in \cite{DemboZeitouni}) and lower
semi-continuous, then $\{Z_n:n\geq 1\}$ satisfies the LDP with
speed function $v_n$ and good rate function $\Lambda^*$ defined by
$\Lambda^*(z):=\sup_{\theta\in\mathbb{R}}\{\theta
z-\Lambda(\theta)\}$. For the sake of completeness we recall that the
function $\Lambda$ is essentially smooth if the interior of 
$\mathcal{D}(\Lambda)$ is non-empty, it is differentiable throughout the
interior of $\mathcal{D}(\Lambda)$, and $|\Lambda^\prime(\theta_n)|\to\infty$
whenever $\{\theta_n\}$ is a sequence of points in the interior of 
$\mathcal{D}(\Lambda)$ which converges to a boundary point of 
$\mathcal{D}(\Lambda)$.

\subsection{Preliminaries on the sequence $\{C_n(w):n\geq 1\}$}\label{sub:model}
Let $\{X_n:n\geq 1\}$ be a sequence of i.i.d. positive random
variables and let $X_{1:n}\leq\cdots\leq X_{n:n}$ be the ascending
order statistics of $X_1,\ldots,X_n$ (for all $n\geq 1$); moreover we
set $X_{0:n}=0$. Then we consider the sequence $\{C_n(w):n\geq 1\}$
defined by
\begin{equation}\label{eq:def-Cnw}
C_n(w):=\sum_{k=0}^{n-1}w(k/n)(X_{k+1:n}-X_{k:n}),
\end{equation}
for some function $w:[0,1]\to\mathbb{R}$. So we have
$$C_n(w)=\sum_{k=0}^{n-1}w(k/n)X_{k+1:n}-\sum_{k=0}^{n-1}w(k/n)X_{k:n}
=\sum_{k=1}^n w((k-1)/n)X_{k:n}-\sum_{k=0}^{n-1}w(k/n)X_{k:n}$$
and, by taking into account $X_{0:n}=0$, we get
\begin{equation}\label{eq:def-without-spacings}
C_n(w)=\sum_{k=1}^{n-1}(w((k-1)/n)-w(k/n))X_{k:n}+w((n-1)/n)X_{n:n}.
\end{equation}

Actually in this paper we assume that the common distribution of
the random variables $\{X_n:n\geq 1\}$ is $\mathcal{EXP}(\lambda)$
for some $\lambda>0$, i.e. their (common) distribution function is
\begin{equation}\label{eq:exponential-df}
F(t):=1-e^{-\lambda t}\ \mbox{for all}\ t\geq 0.
\end{equation}
Then, in such a case, it is known (see e.g. Subsection 2.3 in \cite{Pyke})
that the spacings
$$\{X_{1:n}-X_{0:n},X_{2:n}-X_{1:n},\ldots,X_{n:n}-X_{n-1:n}\}$$
are independent and, for all $k\in\{0,\ldots,n-1\}$, the distribution of
$X_{k+1:n}-X_{k:n}$ is $\mathcal{EXP}(\lambda(n-k))$. This result 
yields some explicit formulas for moment generating function, mean
and variance of $C_n(w)$. Firstly, for all $\theta\in\mathbb{R}$, we have
$$\mathbb{E}\left[e^{\theta C_n(w)}\right]=\prod_{k=0}^{n-1}\mathbb{E}\left[e^{\theta w(k/n)(X_{k+1:n}-X_{k:n})}\right],$$
and therefore
\begin{equation}\label{eq:mgf}
\mathbb{E}\left[e^{\theta C_n(w)}\right]=\left\{\begin{array}{ll}
\prod_{k=0}^{n-1}\frac{\lambda(n-k)}{\lambda(n-k)-\theta w(k/n)}&\
\mbox{if}\ \theta w(k/n)<\lambda(n-k)\ \mbox{for all}\ k\in\{0,\ldots,n-1\}\\
\infty&\ \mbox{otherwise}.
\end{array}\right.
\end{equation}
Moreover
\begin{equation}\label{eq:mean-variance}
\mathbb{E}[C_n(w)]=\frac{1}{\lambda}\sum_{k=0}^{n-1}\frac{w(k/n)}{n-k}\quad
\mbox{and}\quad
\mbox{Var}[C_n(w)]=\frac{1}{\lambda^2}\sum_{k=0}^{n-1}\frac{w^2(k/n)}{(n-k)^2}.
\end{equation}

Now we discuss the almost sure convergence and the asymptotic normality following the lines of some
proofs in \cite{DicrescenzoLongobardiJSPI2009} and \cite{DicrescenzoLongobardiPROCEEDINGS}. We
introduce the following condition.

\begin{condition}\label{condition-DL-approach}
The function $w:[0,1]\to\mathbb{R}$ is continuous and there exist $x_0\in(0,1)$,
$\beta\in(0,1]$ and $c>0$ such that $|w(x)|\leq c(1-x)^\beta$ for all $x\in[1-x_0,1]$.
\end{condition}

We start with a generalization of Proposition 2 in \cite{DicrescenzoLongobardiPROCEEDINGS}.
In view of what follows we recall that Condition \ref{condition-DL-approach} yields
$w(1)=0$, and this condition is needed to have the finiteness of the almost sure limit
$\int_0^\infty w(F(z))dz$ (see \eqref{eq:as-convergence-PROCEEDINGS} below).

\begin{proposition}\label{prop:DL-prop2-generalization}
Assume that Condition \ref{condition-DL-approach} holds.
Let $\{X_n:n\geq 1\}$ be a sequence of i.i.d. positive random variables in $L^p$ for
some $p$ such that $\beta p>1$, with (common) distribution function $F$ possibly different
from the one in \eqref{eq:exponential-df}. Then
\begin{equation}\label{eq:as-convergence-PROCEEDINGS}
C_n(w)\to\int_0^\infty w(F(z))dz\ \mbox{a.s.}\ (\mbox{as}\ n\to\infty).
\end{equation}
\end{proposition}
\begin{proof}
We follow the lines of the proof of Proposition 2 in
\cite{DicrescenzoLongobardiPROCEEDINGS} (see also the proof of Theorem 9 in
\cite{RaoChenVemuriWang}). Obviously we have
$$C_n(w)=\int_0^\infty w(\hat{F}_n(z))dz\ (\mbox{for all}\ n\geq 1),$$
where $\hat{F}_n(x):=\frac{1}{n}\sum_{k=1}^n1_{\{X_k\leq x\}}$ is the empirical distribution
function. We take $a_0>0$ such that $F(a_0)\geq 1-\frac{x_0}{2}$ and, by the Glivenko Cantelli
Theorem, for $n$ large enough we have
$$F(a_0)+\frac{x_0}{2}\geq\hat{F}_n(a_0)\geq F(a_0)-\frac{x_0}{2}.$$
Thus for all $z\geq a_0$ we have
$$\hat{F}_n(z)\geq\hat{F}_n(a_0)\geq 1-x_0,$$
which yields
$$|w(\hat{F}_n(z))|\leq c(1-\hat{F}_n(z))^\beta$$
by Condition \ref{condition-DL-approach}. We also remark that
$$1-\hat{F}_n(z)=1-\frac{1}{n}\sum_{k=1}^n1_{\{X_k\leq z\}}=\frac{1}{n}\sum_{k=1}^n1_{\{X_k>z\}}\leq\frac{1}{n}\sum_{k=1}^n1_{\{X_k>z\}}\frac{X_k^p}{z^p}
\leq\frac{1}{n}\sum_{k=1}^n\frac{X_k^p}{z^p}\leq\frac{\alpha}{z^p},$$
where $\alpha:=\sup_{n\geq 1}\frac{1}{n}\sum_{k=1}^nX_k^p<\infty$ a.s. (in fact, since the
random variables $\{X_n:n\geq 1\}$ are in $L^p$, $\alpha$ is the supremum of a sequence that
converges a.s.); thus
$$|w(\hat{F}_n(z))|\leq c\frac{\alpha^\beta}{z^{\beta p}}.$$
So, by the Glivenko Cantelli Theorem, we can apply the dominated convergence theorem (noting 
that $\int_{a_0}^\infty\frac{dz}{z^{\beta p}}<\infty$ because $\beta p>1$) and we have
$$\int_{a_0}^\infty w(\hat{F}_n(z))dz\to\int_{a_0}^\infty w(F(z))dz\ \mbox{a.s.}\ (\mbox{as}\ n\to\infty).$$
Then we easily conclude the proof noting that we also have
$$\int_0^{a_0}w(\hat{F}_n(z))dz\to\int_0^{a_0}w(F(z))dz\ \mbox{a.s.}\ (\mbox{as}\ n\to\infty)$$
again by the Glivenko Cantelli Theorem and the dominated convergence theorem (noting that
$w$ is continuous and therefore bounded, and the integral is over a bounded interval).
\end{proof}

In particular, if $F$ is the distribution function in \eqref{eq:exponential-df}, it is easy to
check that the limit value is
\begin{equation}\label{eq:limit-expression}
\int_0^\infty w(F(z))dz=\int_0^\infty w(1-e^{-\lambda z})dz=\frac{1}{\lambda}\int_0^1\frac{w(x)}{1-x}dx=:\mu_w,
\end{equation}
which is finite by Condition \ref{condition-DL-approach}; moreover, if we take the
mean value in \eqref{eq:mean-variance}, we have
\begin{equation}\label{eq:limit-mean}
\lim_{n\to\infty}\mathbb{E}[C_n(w)]=\mu_w.
\end{equation}

We conclude with a brief comment on the asymptotic Normality of the empirical
estimators, i.e. the weak convergence of
$\frac{C_n(w)-\mathbb{E}[C_n(w)]}{\sqrt{\mbox{Var}[C_n(w)]}}$ to the standard
Normal distribution. We can follow the lines of the proof of Theorem 7.1 in
\cite{DicrescenzoLongobardiJSPI2009} and, in particular, the Lyapunov condition
for the sequence $\{C_n(w):n\geq 1\}$ is
\begin{equation}\label{eq:Lyapunov-condition}
\lim_{n\to\infty}\frac{\frac{1}{(\lambda n)^3}\sum_{k=0}^{n-1}\frac{|w(k/n)|^3}{(1-k/n)^3}}
{\left(\frac{1}{(\lambda n)^2}\sum_{k=0}^{n-1}\frac{w^2(k/n)}{(1-k/n)^2}\right)^{3/2}}=0.
\end{equation}

\begin{remark}\label{rem:added}
By taking into account Condition \ref{condition-DL-approach}, it is easy to 
check that \eqref{eq:Lyapunov-condition} holds if
\begin{equation}\label{eq:finite-integrals-for-2-and-3}
\lim_{n\to\infty}\frac{1}{n}\sum_{k=0}^{n-1}\frac{|w(k/n)|^3}{(1-k/n)^3}
=\int_0^1\frac{|w(x)|^3}{(1-x)^3}dx<\infty;
\end{equation}
this yields $3(1-\beta)<1$, and therefore $\beta>\frac{2}{3}$.
\end{remark}

\begin{remark}\label{rem:weak-convergence-non-standard-case}
We have
\begin{equation}\label{eq:limit-variance}
\lim_{n\to\infty}n\mbox{Var}[C_n(w)]=\frac{1}{\lambda^2}\int_0^1\frac{w^2(x)}{(1-x)^2}dx=:\sigma_w^2;
\end{equation}
thus the above weak convergence of $\frac{C_n(w)-\mathbb{E}[C_n(w)]}{\sqrt{\mbox{Var}[C_n(w)]}}$
to the standard Normal distribution is equivalent to the weak convergence of
$\sqrt{n}(C_n(w)-\mathbb{E}[C_n(w)])$ to the centered Normal distribution
with variance $\sigma_w^2$.
\end{remark}

Some examples for the function $w$ with $\beta=1$ will be 
presented just after Condition \ref{condition} (see \eqref{eq:illustrative-examples}).
An example with $\beta\in(0,1)$ is $w(x)=(1-x)^\beta$; then, by \eqref{eq:limit-expression}
and \eqref{eq:limit-variance}, we have
$$\mu_w=\frac{1}{\lambda}\int_0^1(1-x)^{\beta-1}dx=\frac{1}{\lambda\beta}$$
and, if $\beta>1/2$,
$$\sigma_w^2=\frac{1}{\lambda^2}\int_0^1(1-x)^{2\beta-2}dx=\frac{1}{\lambda^2(2\beta-1)}.$$

\subsection{Connections with some literature}\label{sub:literature}
We note that the sequence $\{C_n(w):n\geq 1\}$ defined by \eqref{eq:def-Cnw} (see also 
\eqref{eq:def-without-spacings}) coincides with the sequence $\{L_n:n\geq 1\}$ of $L$-statistics
in \cite{GaoZhaoAS2011} (Section 4.6) if we take $w(\cdot)=w(J;\cdot)$, where $w(J;\cdot)$ is
defined by
\begin{equation}\label{eq:def-w-with-score-function}
w(J;x):=\int_x^1 J(u)du
\end{equation}
for some function $J$ called \emph{score function} (a sequence of estimators of this kind appears
in several references; here we recall \cite{JonesZitikis}, eqs. (18) and (19), for the estimation
of risk measures and related quantities). In \cite{GaoZhaoAS2011} it is not required that the 
i.i.d. random variables $\{X_n:n\geq 1\}$ are exponentially distributed.

Moreover, if we consider the score function
$$\tilde{J}(u):=\log u+1,$$
we get
$$w(\tilde{J};x):=\int_x^1 \log u+1du=[u\log u]_{u=x}^{x=1}=-x\log x;$$
then, by \eqref{eq:def-Cnw} (and by taking into account that $0\log 0=0$),
we get
$$C_n(w(\tilde{J};\cdot))=\sum_{k=1}^{n-1}\left(-\frac{k}{n}\log \frac{k}{n}\right)(X_{k+1:n}-X_{k:n}).$$
So $\{C_n(w(\tilde{J};\cdot)):n\geq 1\}$ coincides with:
\begin{itemize}
\item $\{\mathcal{CE}(\hat{F}_n):n\geq 1\}$ in \cite{DicrescenzoLongobardiJSPI2009} (Section 7),
when $\{X_n:n\geq 1\}$ are i.i.d. and positive random variables;
\item $\{C_n:n\geq 1\}$ in \cite{GiulianoMacciPacchiarottiJSPI2015} (Section 4), when
$\{X_n:n\geq 1\}$ are i.i.d. $\mathcal{EXP}(\lambda)$ distributed random variables.
\end{itemize}

\section{Results}\label{sec:results}
In this section we generalize the results for the sequence $\{C_n:n\geq 1\}$ in
\cite{GiulianoMacciPacchiarottiJSPI2015} (Section 4). In view of what follows we
introduce the following condition.

\begin{condition}\label{condition}
Let $w:[0,1]\to\mathbb{R}$ be a function as in Condition \ref{condition-DL-approach}
with $\beta=1$, and set $h_w(x):=\frac{w(x)}{1-x}$ for $x\in[0,1)$. Moreover let
$\Lambda_w:\mathbb{R}\to\mathbb{R}\cup\{\infty\}$ be the function defined by
$$\Lambda_w(\theta):=\left\{\begin{array}{ll}
\int_0^1\log\left(\frac{\lambda}{\lambda-\theta h_w(x)}\right)dx&\ \mathrm{if}\ \sup_{x\in[0,1)}\{\theta h_w(x)\}\leq\lambda\\
\infty&\ \mathrm{otherwise},
\end{array}\right.$$
and assume that $\Lambda_w$ is finite in a neighbourhood of the origin $\theta=0$.
\end{condition}

We remark that the function $\Lambda_w$ would not be finite in a neighbourhood of the
origin $\theta=0$ if we have Condition \ref{condition-DL-approach} with $\beta\in(0,1)$.

\paragraph{Some examples for the function $w$.}
We consider the following functions:
\begin{equation}\label{eq:illustrative-examples}
w_1(x):=1-x;\quad w_2(x):=(1-x)^2;\quad w_3(x):=(1-x)(1-\sqrt{x}).
\end{equation}
For all these cases Condition \ref{condition-DL-approach} holds with $\beta=1$;
moreover: $\sup_{x\in[0,1)}\{\theta h_{w_1}(x)\}\leq\lambda$ if and only if $\theta\leq\lambda$,
$\Lambda_w$ is lower semicontinuous, and there exists $\Lambda_w^\prime(\theta)$ for
$\theta<\lambda$. Thus, for each function, we have to check the steepness of $\Lambda_w$,
i.e.
\begin{equation}\label{eq:steepness-condition-for-examples}
\lim_{\theta\to\lambda^-}\Lambda_w^\prime(\theta)=\infty,
\end{equation}
which yields its essential smoothness required in the statement of Proposition \ref{prop:LDP}.
\begin{itemize}
\item For $w=w_1$ we have
\begin{equation}\label{eq:Lambda-for-the-case-of-Cramer-theorem}
\Lambda_{w_1}(\theta)=\left\{\begin{array}{ll}
\log\left(\frac{\lambda}{\lambda-\theta}\right)&\ \mathrm{if}\ \theta<\lambda\\
\infty&\ \mathrm{otherwise}.
\end{array}\right.
\end{equation}
So we have $\Lambda_{w_1}(\lambda)=\infty$, and therefore
\eqref{eq:steepness-condition-for-examples} holds; indeed we have
$$\lim_{\theta\to\lambda^-}\Lambda_{w_1}^\prime(\theta)=\lim_{\theta\to\lambda^-}\frac{1}{\lambda-\theta}=\infty.$$
\item For $w=w_2$ we have
$$\Lambda_{w_2}(\lambda)=-\int_0^1\log(1-h_{w_2}(x))dx=1;$$
however, even if $\Lambda_{w_2}(\lambda)<\infty$, \eqref{eq:steepness-condition-for-examples}
holds because
$$\lim_{\theta\to\lambda^-}\Lambda_{w_2}^\prime(\theta)=\lim_{\theta\to\lambda^-}\int_0^1\frac{h_{w_2}(x)}{\lambda-\theta h_{w_2}(x)}dx
=\frac{1}{\lambda}\int_0^1\frac{h_{w_2}(x)}{1-h_{w_2}(x)}dx=\frac{1}{\lambda}\left(\int_0^1\frac{1}{x}dx-1\right)=\infty.$$
\item For $w=w_3$ we have
$$\Lambda_{w_3}(\lambda)=-\int_0^1\log(1-h_{w_3}(x))dx=\frac{1}{2}$$
and
$$\lim_{\theta\to\lambda^-}\Lambda_{w_3}^\prime(\theta)=\lim_{\theta\to\lambda^-}\int_0^1\frac{h_{w_3}(x)}{\lambda-\theta h_{w_3}(x)}dx
=\frac{1}{\lambda}\int_0^1\frac{h_{w_3}(x)}{1-h_{w_3}(x)}dx=\frac{1}{\lambda}\left(\int_0^1\frac{1}{\sqrt{x}}dx-1\right)=\frac{1}{\lambda};$$
thus \eqref{eq:steepness-condition-for-examples} fails.
\end{itemize}

We start with the first result, which is the analogue of Proposition
4.1 in \cite{GiulianoMacciPacchiarottiJSPI2015}.

\begin{proposition}\label{prop:LDP}
Assume that $\{X_n:n\geq 1\}$ are i.i.d. and $\mathcal{EXP}(\lambda)$
distributed, Condition \ref{condition} holds, and $\Lambda_w$ is
essentially smooth and lower semi-continuous. Then the sequence
$\{C_n(w):n\geq 1\}$ defined by \eqref{eq:def-Cnw} satisfies
the LDP with speed function $v_n=n$ and good rate function $\Lambda_w^*$
defined by
$$\Lambda_w^*(y):=\sup_{\theta\in\mathbb{R}}\{\theta y-\Lambda_w(\theta)\}.$$
\end{proposition}
\begin{proof}
We want to apply G\"{a}rtner Ellis Theorem; thus we have to check that
\begin{equation}\label{eq:GElimitLD}
\lim_{n\to\infty}\frac{1}{n}\log\mathbb{E}\left[e^{n\theta
C_n(w)}\right]=\Lambda_w(\theta)\ (\mbox{for all}\ \theta\in\mathbb{R}).
\end{equation}
We remark that, by \eqref{eq:mgf}, we have
$$\frac{1}{n}\log\mathbb{E}\left[e^{n\theta C_n(w)}\right]
=\frac{1}{n}\sum_{k=0}^{n-1}\log\left(\frac{\lambda(n-k)}{\lambda(n-k)-n\theta w(k/n)}\right)
=\frac{1}{n}\sum_{k=0}^{n-1}\log\left(\frac{\lambda\left(1-\frac{k}{n}\right)}{\lambda\left(1-\frac{k}{n}\right)-\theta w(k/n)}\right)$$
for all $\theta\in\mathbb{R}$ such that
\begin{equation}\label{eq:constraint}
\theta w(k/n)<\lambda\left(1-\frac{k}{n}\right)\
\mbox{for all}\ j\in\{0,\ldots,n-1\}
\end{equation}
(and $\frac{1}{n}\log\mathbb{E}\left[e^{n\theta C_n(w)}\right]$ equal
to infinity otherwise). Moreover condition \eqref{eq:constraint} holds
(for any fixed $n\geq 1$) if and only if
$$\theta h_w(k/n)<\lambda\ \mbox{for all}\ k\in\{0,\ldots,n-1\}.$$
Thus the limit in \eqref{eq:GElimitLD} trivially holds if
$\sup_{x\in[0,1)}\{\theta h_w(x)\}>\lambda$ while, if
$\sup_{x\in[0,1)}\{\theta h_w(x)\}\leq\lambda$, the limit
\eqref{eq:GElimitLD} can be checked noting that we have a limit of an
integral sum (possibly equal to infinity). In conclusion the desired
LDP holds as a straightforward application of the G\"{a}rtner Ellis
Theorem.
\end{proof}

\begin{remark}\label{rem:unique-zero}
It is well-known that $\Lambda_w^*(y)=0$ if and only if
$y=\Lambda_w^\prime(0)$. Then, since we can differentiate under
the integral sign by Condition \ref{condition}, we get
$$\Lambda_w^\prime(0)=\frac{1}{\lambda}\int_0^1h_w(x)dx=\frac{1}{\lambda}\int_0^1\frac{w(x)}{1-x}dx,$$
i.e. $\Lambda_w^\prime(0)$ coincides with $\mu_w$ in \eqref{eq:limit-expression}.
\end{remark}

\begin{remark}\label{rem:connections-with-najim-paper}
If the random variables in Proposition \ref{prop:LDP} are not exponentially distributed, then we cannot
rely on some properties of the spacings cited above (independence and exponential distributions with 
different parameters); so we have some difficulties to apply the G\"{a}rtner Ellis Theorem. A possible 
way to overcome this problem is to try to apply Theorem 2.2 in \cite{Najim}. Some technical conditions
should be checked and this could be done in a successive work.
\end{remark}

\begin{remark}\label{rem:connections-with-cramer-theorem}
Here we consider Proposition \ref{prop:LDP} with $w=w_1$, where $w_1$ is the
function in \eqref{eq:illustrative-examples}. Thus $\Lambda_w$ coincides with
the function $\Lambda_{w_1}$ in \eqref{eq:Lambda-for-the-case-of-Cramer-theorem};
moreover we can check (after some easy computations) that $\Lambda_w^*$ coincides
with
$$\Lambda_{w_1}^*(y)=\left\{\begin{array}{ll}
\lambda y-1-\log(\lambda y)&\ \mbox{if}\ y>0\\
\infty&\ \mbox{otherwise}.
\end{array}\right.$$
Then we have the rate function provided by the Cram\'er Theorem (see e.g.
Theorem 2.2.3 in \cite{DemboZeitouni}) for the sequence of empirical means
$\left\{\frac{X_1+\cdots+X_n}{n}:n\geq 1\right\}$ when (as happens in
Proposition \ref{prop:LDP}) $\{X_n:n\geq 1\}$ is a sequence of i.i.d. and
$\mathcal{EXP}(\lambda)$ distributed random variables. In fact it is easy to
check that
$$C_n(w_1)=\frac{1}{n}\sum_{k=0}^{n-1}(n-k)(X_{k+1:n}-X_{k:n})\ (\mbox{for all}\ n\geq 1)$$
by \eqref{eq:def-Cnw} and the definition of $w_1$ in
\eqref{eq:illustrative-examples}, and therefore $\left\{C_n(w_1):n\geq 1\right\}$
and $\left\{\frac{X_1+\cdots+X_n}{n}:n\geq 1\right\}$ are equally distributed
by taking into account the independence and the distributions of the spacings
(indeed, for each $n\geq 1$, the law of $C_n(w_1)$ and
$\frac{X_1+\cdots+X_n}{n}$ is the Gamma distribution with probability density
function $g(z)=\frac{\lambda^n}{(n-1)!}z^{n-1}e^{-\lambda z}1_{(0,\infty)}(z)$).
\end{remark}

The second result, which is the analogue of Proposition 4.2 in
\cite{GiulianoMacciPacchiarottiJSPI2015}, provides an upper bound of the rate
function $\Lambda_w^*$ in Proposition \ref{prop:LDP} when $h_w(x)>0$
almost everywhere with respect to $x$. This upper bound can be expressed in
terms of the \emph{relative entropy} (see e.g. \cite{KL1951}) of an exponential
distribution with respect to another one. We recall that, given two absolutely
continuous real valued random variables $X_1$ and $X_2$ with densities $f_1$ and
$f_2$, the relative entropy of $X_1$ with respect to $X_2$ is defined by
$$H(X_1|X_2):=\int_{\mathbb{R}}f_1(x)\log\frac{f_1(x)}{f_2(x)}dx;$$
thus $H(X_1|X_2)$ actually depends on the laws of the random variables $X_1$
and $X_2$. Then the relative entropy of the distribution $\mathcal{EXP}(\lambda_1)$
with respect to the distribution $\mathcal{EXP}(\lambda_2)$ is
$$H(\mathcal{EXP}(\lambda_1)|\mathcal{EXP}(\lambda_2))=\frac{\lambda_2}{\lambda_1}-1-\log\frac{\lambda_2}{\lambda_1}.$$

\begin{proposition}\label{prop:on-rf-in-LDP}
Let $h_w$ be as in Condition \ref{condition} and assume that
$h_w(x)>0$ almost everywhere with respect to $x$. Moreover set
$M_w(y):=\int_0^1H\left(\mathcal{EXP}(1/y)|\mathcal{EXP}(\lambda
h_w^{-1}(x)\right)dx$. Then: (i) $\Lambda_w^*(y)\leq M_w(y)$ for all
$y\in(0,\infty)$; (ii) $\Lambda_w^*(y)=\infty$ for all
$y\in (-\infty,0]$; (iii) the infimum of $M_w(y)$ is attained at
$y=\bar{y}_w$, where $\bar{y}_w:=(\lambda\int_0^1h_w^{-1}(x)dx)^{-1}$.
\end{proposition}
\begin{proof}
We start with the proof of (i). We remark that, for
$y>0$, we have
$$\sup_{\theta<\eta}\left\{\theta y-\log\left(\frac{\eta}{\eta-\theta}\right)\right\}=H(\mathcal{EXP}(1/y)|\mathcal{EXP}(\eta))$$
for $\eta>0$; then we get
\begin{eqnarray*}
\Lambda_w^*(y)&=&\sup_{\theta\sup_{z\in [0,1)}h_w(z)\leq\lambda}\left\{\theta y-
\int_0^1\log\left(\frac{\lambda}{\lambda-\theta h_w(x)}\right)dx\right\}\\
&=&\sup_{\theta\sup_{z\in [0,1)}h_w(z)\leq\lambda}\left\{\theta y-
\int_0^1\log\left(\frac{\lambda h_w^{-1}(x)}{\lambda h_w^{-1}(x)-\theta}\right)dx\right\}\\
&\leq&\int_0^1\sup_{\theta\sup_{z\in [0,1)}h_w(z)\leq\lambda}\left\{\theta
y-\log\left(\frac{\lambda h_w^{-1}(x)}{\lambda h_w^{-1}(x)-\theta}\right)\right\}dx\\
&\leq&\int_0^1\sup_{\theta<\lambda h_w^{-1}(x)}\left\{\theta
y-\log\left(\frac{\lambda h_w^{-1}(x)}{\lambda h_w^{-1}(x)-\theta}\right)\right\}dx
=\int_0^1H\left(\mathcal{EXP}(1/y)|\mathcal{EXP}(\lambda h_w^{-1}(x))\right)dx.
\end{eqnarray*}
Now the proof of (ii): for $y<0$ we have
$$\Lambda_w^*(y)\geq\sup_{\theta\leq 0}\left\{\theta y-\int_0^1\log\left(\frac{\lambda h_w^{-1}(x)}{\lambda h_w^{-1}(x)-\theta}\right)dx\right\}
\geq\sup_{\theta\leq 0}\{\theta y\}=\infty;$$
for $y=0$ (this case was forgotten in the proof of Proposition 4.2
in \cite{GiulianoMacciPacchiarottiJSPI2015}) we have
$$\Lambda_w^*(0)\geq\sup_{\theta\leq 0}\left\{-\int_0^1\log\left(\frac{\lambda h_w^{-1}(x)}{\lambda h_w^{-1}(x)-\theta}\right)dx\right\}
=\lim_{\theta\to-\infty}-\int_0^1\log\left(\frac{\lambda h_w^{-1}(x)}{\lambda h_w^{-1}(x)-\theta}\right)dx=\infty.$$
Finally the proof of (iii). One can check that
$$M_w(y)=\lambda\int_0^1h_w^{-1}(x)dx\cdot y-1-\log\lambda-\int_0^1\log(h_w^{-1}(x))dx-\log y$$
and its derivative is
$$M_w^\prime(y)=\lambda\int_0^1h_w^{-1}(x)dx-\frac{1}{y}.$$
So we have $M_w^\prime(y)=0$ if and only if $y=\bar{y}_w$, and
$y=\bar{y}_w$ is a minimizer by the convexity of $M_w$.
\end{proof}

The third result, which is the analogue of Proposition 4.3 in
\cite{GiulianoMacciPacchiarottiJSPI2015}, concerns moderate
deviations. In view of its proof we remark that
\begin{equation}\label{eq:local-estimate-degree3}
\mbox{there exists}\ \delta>0\ \mbox{such that}\ \log(1+x)\leq
x-\frac{x^2}{2}+\frac{x^3}{3}\ \mbox{for all}\ |x|<\delta
\end{equation}
(which can be proved by checking that the function $g$ defined by
$g(x):=\log(1+x)-(x-\frac{x^2}{2}+\frac{x^3}{3})$ has a local
maximum at $x=0$) and
\begin{equation}\label{eq:local-estimate-degree2}
\mbox{for all}\ v>\frac{1}{2},\ \mbox{there exists}\ \delta>0\
\mbox{such that}\ \log(1+x)\geq x-vx^2\ \mbox{for all}\ |x|<\delta
\end{equation}
(which can be proved by checking that the function $g$ defined by
$g(x):=\log(1+x)-(x-vx^2)$ has a local minimum at $x=0$).

\begin{proposition}\label{prop:MD}
Assume that $\{X_n:n\geq 1\}$ are i.i.d. and $\mathcal{EXP}(\lambda)$
distributed, and Condition \ref{condition} holds. Then, for any
positive sequence $\{a_n:n\geq 1\}$ such that
\begin{equation}\label{eq:MD-conditions}
a_n\to 0\quad\mbox{and}\quad na_n\to\infty\ (\mbox{as}\ n\to\infty),
\end{equation}
the sequence $\left\{\sqrt{na_n}(C_n(w)-\mathbb{E}[C_n(w)]):n\geq 1\right\}$
satisfies the LDP with speed function $v_n=1/a_n$ and good rate
function $\tilde{\Lambda}_w^*(y)$ defined by
$\tilde{\Lambda}_w^*(y):=\frac{y^2}{2\sigma_w^2}$, where
$\sigma_w^2$ is the expression in \eqref{eq:limit-variance}.
\end{proposition}
\begin{proof}
We want to apply the G\"{a}rtner Ellis Theorem with speed function $1/a_n$; thus
we have to check that
\begin{equation}\label{eq:GElimitMD-LB}
\liminf_{n\to\infty}a_n\log\mathbb{E}\left[\exp\left(\theta\sqrt{\frac{n}{a_n}}
(C_n(w)-\mathbb{E}[C_n(w)])\right)\right]\geq\sigma_w^2\frac{\theta^2}{2}
\end{equation}
and
\begin{equation}\label{eq:GElimitMD-UB}
\limsup_{n\to\infty}a_n\log\mathbb{E}\left[\exp\left(\theta\sqrt{\frac{n}{a_n}}
(C_n(w)-\mathbb{E}[C_n(w)])\right)\right]\leq\sigma_w^2\frac{\theta^2}{2}
\end{equation}
for all $\theta\in\mathbb{R}$.

It is useful to remark that, by \eqref{eq:mgf} and the mean value in
\eqref{eq:mean-variance} (together with some computations), we have
\begin{multline*}
\log\mathbb{E}\left[\exp\left(\theta\sqrt{\frac{n}{a_n}}(C_n(w)-\mathbb{E}[C_n(w)])\right)\right]
=\log\mathbb{E}\left[e^{\theta\sqrt{\frac{n}{a_n}}C_n(w)}\right]
-\theta\sqrt{\frac{n}{a_n}}\mathbb{E}[C_n(w)]\\
=\sum_{k=0}^{n-1}\log\frac{\lambda(n-k)}{\lambda(n-k)-\theta\sqrt{\frac{n}{a_n}}w(k/n)}
-\frac{\theta}{\lambda}\sqrt{\frac{n}{a_n}}\sum_{k=0}^{n-1}\frac{w(k/n)}{n-k}\\
=-\sum_{k=0}^{n-1}\left(\log\left(1-\frac{\theta}{\lambda\sqrt{na_n}}\frac{w(k/n)}{1-k/n}\right)
+\frac{\theta}{\lambda\sqrt{na_n}}\frac{w(k/n)}{1-k/n}\right)
\end{multline*}
for all $\theta\in\mathbb{R}$ such that
$$\frac{\theta}{\lambda\sqrt{na_n}}\frac{w(k/n)}{1-k/n}<1\ \mbox{for all}\ k\in\{0,\ldots,n-1\}$$
(and $\log\mathbb{E}\left[\exp\left(\theta\sqrt{\frac{n}{a_n}}(C_n(w)-\mathbb{E}[C_n(w)])\right)\right]$
equal to infinity otherwise). Then, by Condition \ref{condition} and
by $na_n\to\infty$, for all $\delta>0$ there exists $\bar{n}$ such that
$$\left|\frac{\theta}{\lambda\sqrt{na_n}}\frac{w(k/n)}{1-k/n}\right|<\delta\ \mbox{for all}\ k\in\{0,\ldots,n-1\}$$
for all $n>\bar{n}$ (in fact $\left|\frac{\theta}{\lambda\sqrt{na_n}}\frac{w(k/n)}{1-k/n}\right|
\leq\frac{|\theta|c}{\lambda\sqrt{na_n}}\to 0$ as $n\to\infty$).

Now we are ready for the proof of \eqref{eq:GElimitMD-LB} and
\eqref{eq:GElimitMD-UB}; this will be done by using
\eqref{eq:local-estimate-degree3} and \eqref{eq:local-estimate-degree2}
for $\delta>0$ chosen above and for suitable choices of $x$ which depend on
$n>\bar{n}$. We start with the proof of \eqref{eq:GElimitMD-LB}. If we combine
the above computations in this proof and \eqref{eq:local-estimate-degree3} (with
$x=-\frac{\theta}{\lambda\sqrt{na_n}}\frac{w(k/n)}{1-k/n}$), we have
$$a_n\log\mathbb{E}\Big[\exp\Big(\theta\sqrt{\frac{n}{a_n}}(C_n(w)-\mathbb{E}[C_n(w)])\Big)\Big]
\geq a_n\sum_{k=0}^{n-1}\Big(\frac{\theta^2}{2\lambda^2na_n}\frac{w^2(k/n)}{(1-k/n)^2}
+\frac{\theta^3}{3\lambda^3(na_n)^{3/2}}\frac{w^3(k/n)}{(1-k/n)^3}\Big);$$
hence, by taking into account the limit for the variance in
\eqref{eq:limit-variance} and
$$\lim_{n\to\infty}\frac{1}{n\sqrt{na_n}}\sum_{k=0}^{n-1}\frac{w^3(k/n)}{(1-k/n)^3}=0$$
(because $na_n\to\infty$ by \eqref{eq:MD-conditions} and, as explained
in Remark \ref{rem:added},
$\lim_{n\to\infty}\frac{1}{n}\sum_{k=0}^{n-1}\frac{w^3(k/n)}{(1-k/n)^3}=\int_0^1\frac{w^3(x)}{(1-x)^3}dx$
is finite because $\beta>\frac{2}{3}$), we obtain
$$\liminf_{n\to\infty}a_n\log\mathbb{E}\left[\exp\left(\theta\sqrt{\frac{n}{a_n}}(C_n(w)-\mathbb{E}[C_n(w)])\right)\right]
\geq\sigma_w^2\frac{\theta^2}{2}$$ and
\eqref{eq:GElimitMD-LB} is proved. The proof of \eqref{eq:GElimitMD-UB} is
similar. We have to consider \eqref{eq:local-estimate-degree2} instead of
\eqref{eq:local-estimate-degree3} and, again, we take into account the limit
of the variance \eqref{eq:limit-variance}; then we obtain
$$\limsup_{n\to\infty}a_n\log\mathbb{E}\left[\exp\left(\theta\sqrt{\frac{n}{a_n}}(C_n(w)-\mathbb{E}[C_n(w)])\right)\right]
\leq\limsup_{n\to\infty}a_n\sum_{k=0}^{n-1}\frac{v\theta^2}{na_n}\frac{w^2(k/n)}{(1-k/n)^2}=\sigma_w^2v\theta^2,$$
and we get \eqref{eq:GElimitMD-UB} by letting $v$ go to
$\frac{1}{2}$.
\end{proof}

In the following remark we recall some typical features on moderate deviations.

\begin{remark}\label{rem:comments-collection}
The class of LDPs in Proposition \ref{prop:MD} fill the gap
between two asymptotic regimes.
\begin{enumerate}
\item The almost sure convergence of $C_n(w)$ to $\mu_w$, which is equivalent (by
\eqref{eq:limit-mean}) to the almost sure convergence of $C_n(w)-\mathbb{E}[C_n(w)]$
to zero.
\item The weak convergence of $\sqrt{n}(C_n(w)-\mathbb{E}[C_n(w)])$
to the centered Normal distribution with variance $\sigma_w^2$ (see
Remark \ref{rem:weak-convergence-non-standard-case}).
\end{enumerate}
Then we recover these two cases by taking the sequence of random
variables in Proposition \ref{prop:MD} with $a_n=\frac{1}{n}$ and
$a_n=1$, respectively; in both cases one condition in
\eqref{eq:MD-conditions} holds, and the other one fails.

Moreover we know that the LDP in Proposition \ref{prop:LDP}, which concerns the almost
sure convergence of $C_n(w)$ to $\mu_w$, is governed by the rate function $\Lambda_w^*(y)$
which uniquely vanishes at $y=\Lambda_w^\prime(0)$ (see Remark \ref{rem:unique-zero}), and
$(\Lambda_w^*)^{\prime\prime}(\Lambda_w^\prime(0))=(\Lambda_w^{\prime\prime}(0))^{-1}$. So,
since we can differentiate (twice) under the integral sign by Condition \ref{condition}, we
get (see also \eqref{eq:limit-variance})
$$\Lambda_w^{\prime\prime}(0)=\frac{1}{\lambda}\int_0^1h_w^2(x)dx=\frac{1}{\lambda^2}\int_0^1\frac{w^2(x)}{(1-x)^2}dx=\sigma_w^2,$$
i.e. the variance of the weak limit law of
$\frac{C_n(w)-\mathbb{E}[C_n(w)]}{\sqrt{n}}$.

In some sense we can say that we have an asymptotic normality result
as a consequence of an LDP; an interesting discussion on this issue
can be found in \cite{BrycSPL1993}.
\end{remark}

Finally we show how to obtain a lower bound for the asymptotic variance
$\sigma_w^2$ in Remark \ref{rem:comments-collection} (and in Remark
\ref{rem:weak-convergence-non-standard-case}).

\begin{remark}\label{rem:lower-bound-for-asymptotic-variance}
Here we assume that $\gamma_w:=\int_0^1\frac{w(x)}{1-x}dx\neq 0$.
Then, by \eqref{eq:limit-variance} and an easy application of the
Jensen's inequality, we have
$$\sigma_w^2\geq\frac{1}{\lambda^2}\left(\int_0^1\frac{w(x)}{1-x}dx\right)^2=\frac{\gamma_w^2}{\lambda^2}.$$
So, if we consider the function $w_1$ in \eqref{eq:illustrative-examples},
the inequality turns into an equality if and only if
$$w(x)=\gamma_ww_1(x)=\gamma_w(1-x).$$
From now on we set $\ell_w(x):=\gamma_w(1-x)$; moreover we take
$\gamma_w>0$ and we follow the same lines of some parts of Remark
\ref{rem:connections-with-cramer-theorem}.  Firstly we have
$\Lambda_{\ell_w}(\theta)=\Lambda_{w_1}(\theta\gamma_w)$
for all $\theta\in\mathbb{R}$ and
\begin{multline*}
\Lambda_{\ell_w}^*(y)=\sup_{\theta\in\mathbb{R}}\{\theta y-\Lambda_{\ell_w}(\theta)\}
=\sup_{\theta\in\mathbb{R}}\{\theta y-\Lambda_{w_1}(\theta\gamma_w)\}\\
=\Lambda_{w_1}^*(y\gamma_w^{-1})=\left\{\begin{array}{ll}
\lambda\gamma_w^{-1}y-1-\log(\lambda\gamma_w^{-1}y)&\ \mbox{if}\ y>0\\
\infty&\ \mbox{otherwise}.
\end{array}\right.
\end{multline*}
Moreover $\Lambda_{\ell_w}^*$ is the rate function provided by the Cram\'er Theorem
for a sequence of empirical means of i.i.d. and $\mathcal{EXP}(\lambda\gamma_w^{-1})$
distributed random variables; indeed we have
$$C_n(\ell_w)=\frac{\gamma_w}{n}\sum_{k=0}^{n-1}(n-k)(X_{k+1:n}-X_{k:n})\ (\mbox{for all}\ n\geq 1)$$
by \eqref{eq:def-Cnw} and the definition of $\ell_w$, and therefore
$\left\{C_n(\ell_w):n\geq 1\right\}$ is a sequence of such empirical means
by the independence and the distributions of the spacings.
\end{remark}

\section{Some analogies with a moderate deviation result for $L$-statistics}\label{sec:connections-with-Gao-Zhao}
In this section we discuss some connections between Theorem 4.8 in \cite{GaoZhaoAS2011}
with $F$ as the exponential distribution function in \eqref{eq:exponential-df}, and
Proposition \ref{prop:MD} in this paper with $w(\cdot)=w(J;\cdot)$ as in
\eqref{eq:def-w-with-score-function}.

Firstly, since $F$ is the exponential distribution function in
\eqref{eq:exponential-df}, we can give the following formulas for $m(J,F)$ and
$\sigma^2(J,F)$ in Theorem 4.8 in \cite{GaoZhaoAS2011}:
$$m(J,F):=\int_0^\infty xJ(1-e^{-\lambda x})\lambda e^{-\lambda x}dx;$$
$$\sigma^2(J,F):=\int_0^\infty\int_0^\infty J(1-e^{-\lambda x})J(1-e^{-\lambda y})
\{1-e^{-\lambda(x\wedge y)}-(1-e^{-\lambda x})(1-e^{-\lambda y})\}dxdy.$$
Then, under suitable hypotheses (some of them concern the score function $J$),
Theorem 4.8 in \cite{GaoZhaoAS2011} allows to say that, for any sequence
$\{a(n):n\geq 1\}$ of positive numbers such that
\begin{equation}\label{eq:MD-conditions-GZ}
a(n)\to \infty\quad\mbox{and}\quad \frac{a(n)}{\sqrt{n}}\to 0\ (\mbox{as}\ n\to\infty),
\end{equation}
the sequence $\left\{\frac{\sqrt{n}}{a(n)}(C_n(w)-m(J,F)):n\geq 1\right\}$ satisfies
the LDP with speed $(a(n))^2$ and good rate function
$$I^L(y):=\frac{y^2}{2\sigma^2(J,F)}.$$
Thus $(a(n))^{-2}$ in \eqref{eq:MD-conditions-GZ} plays the role of $a_n$ in
\eqref{eq:MD-conditions}; moreover, as typically happens for the results on moderate
deviations, both rate functions $\Lambda_w^*$ in Proposition \ref{prop:MD} and $I^L$
are quadratic functions that uniquely vanish at the origin $y=0$.

We remark that, if we compare $\left\{\frac{\sqrt{n}}{a(n)}(C_n(w)-m(J,F)):n\geq 1\right\}$
and the sequence of random variables in Proposition \ref{prop:MD} in this paper, by taking
into account the limit \eqref{eq:limit-mean} we expect to have $m(J,F)=\mu_{w(J;\cdot)}$.
In fact, by considering the change of variable $r=1-e^{-\lambda x}$ and some computations
with an integration by parts, we have
\begin{eqnarray*}
m(J,F)&=&\int_0^1-\frac{\log(1-r)}{\lambda}J(r)\lambda(1-r)\frac{dr}{\lambda(1-r)}
=-\frac{1}{\lambda}\int_0^1\log(1-r)J(r)dr\\
&=&-\frac{1}{\lambda}\left\{\left[-w(J;r)\log(1-r)\right]_{r=0}^{r=1}-\int_0^1\frac{w(J;r)}{1-r}dr\right\}
=\frac{1}{\lambda}\int_0^1\frac{w(J;r)}{1-r}dr=\mu_{w(J;\cdot)},
\end{eqnarray*}
indeed $\left[-w(J;r)\log(1-r)\right]_{r=0}^{r=1}=0$ because the score function $J$ is
bounded, continuous and trimmed (i.e. it is equal to zero near $r=0$ and $r=1$).

We also remark that, if we compare the rate functions $\Lambda_w^*$ and $I^L$, we expect
to have $\sigma^2(J,F)=\sigma_{w(J;\cdot)}^2$. In order to check this equality we 
note that the function inside the integral is symmetric with respect to $(x,y)$;
therefore we have the integral over $\{(x,y):0\leq x\leq y\}$ multiplied by 2 and, after
some computations, we get
$$\sigma^2(J,F)=2\int_0^\infty dyJ(1-e^{-\lambda y})e^{-\lambda y}\int_0^y dx (1-e^{-\lambda x})J(1-e^{-\lambda x}).$$
Moreover we consider two further changes of variables: the first one is $r=1-e^{-\lambda x}$,
and we obtain
$$\sigma^2(J,F)=2\int_0^\infty dye^{-\lambda y}J(1-e^{-\lambda y})\int_0^{1-e^{-\lambda y}}\frac{dr}{\lambda(1-r)}rJ(r);$$
the second one is $s=1-e^{-\lambda y}$, and we get
$$\sigma^2(J,F)=2\int_0^1\frac{ds}{\lambda(1-s)}(1-s)J(s)\int_0^s\frac{dr}{\lambda(1-r)}rJ(r)
=\frac{2}{\lambda^2}\int_0^1dsJ(s)\int_0^sdr\frac{r}{1-r}J(r).$$
Finally we conclude with the following computations (again with integration by
parts):
\begin{eqnarray*}
\sigma^2(J,F)&=&\frac{2}{\lambda^2}\left\{\left[-w(J;s)\int_0^sdr\frac{r}{1-r}J(r)\right]_{s=0}^{s=1}+
\int_0^1w(J;s)J(s)\frac{s}{1-s}ds\right\}\\
&=&\frac{2}{\lambda^2}\left\{\left[-\frac{(w(J;s))^2}{2}\frac{s}{1-s}\right]_{s=0}^{s=1}+
\int_0^1\frac{(w(J;s))^2}{2}\frac{1}{(1-s)^2}ds\right\}\\&=&
\frac{1}{\lambda^2}\int_0^1\frac{(w(J;s))^2}{(1-s)^2}ds=\sigma_{w(J;\cdot)}^2,
\end{eqnarray*}
because $\left[-w(J;s)\int_0^sdr\frac{r}{1-r}J(r)\right]_{s=0}^{s=1}=0$ and
$\left[-\frac{(w(J;s))^2}{2}\frac{s}{1-s}\right]_{s=0}^{s=1}=0$ by the hypotheses on
the score function $J$ recalled above and, for the second equality, by Condition
\ref{condition} with $\beta=1$ (which refers to Condition \ref{condition-DL-approach})
for $w(J;\cdot)$.

\begin{remark}\label{rem:MD-for-not-exponential-case}
If the random variables in Proposition \ref{prop:MD} are not exponentially distributed, we have some
difficulties to apply the the G\"{a}rtner Ellis Theorem (as we said in Remark 
\ref{rem:connections-with-najim-paper} for Proposition \ref{prop:LDP}). However one could try to 
adapt the proof of Theorem 4.8 in \cite{GaoZhaoAS2011} (which is proved for $w(\cdot)=w(J;\cdot)$ as in
\eqref{eq:def-w-with-score-function} under some suitable hypotheses) for a quite general function $w$.
This could be done in a successive work.
\end{remark}

\section{Applications to some empirical entropies}\label{sec:empirical-estimators}
A natural way to estimate a functional $\varphi(F)$ of a distribution function $F$ is to consider
$\varphi(\hat{F}_n)$ where, given a sequence $\{X_n:n\geq 1\}$ of i.i.d. random variables with
distribution function $F$ (possibly different from that one in \eqref{eq:exponential-df}),
$\{\hat{F}_n:n\geq 1\}$ is the sequence of the empirical distribution functions defined by
$$\hat{F}_n(x):=\frac{1}{n}\sum_{k=1}^n1_{\{X_k\leq x\}},\quad x \in \mathbb{R}.$$
In this section we concentrate our attention on functionals related to the concept of entropy and
some other related items.

We recall some preliminaries and we refer to \cite{DicrescenzoLongobardiJSPI2009} (see also the
references cited therein). The \emph{cumulative entropy} associated to an absolutely continuous
distribution function $F$ is defined by
$$\mathcal{CE}(F)=-\int_0^\infty F(z)\log F(z)dz.$$
Then, given a sequence of i.i.d. random variables $\{X_n:n\geq 1\}$ with (common) distribution
function $F$, we can consider the sequence of empirical cumulative entropies
$\{\mathcal{CE}(\hat{F}_n):n\geq 1\}$ defined by
$$\mathcal{CE}(\hat{F}_n):=-\int_0^\infty\hat{F}_n(z)\log\hat{F}_n(z)dz.$$
It is known that $\mathcal{CE}(\hat{F}_n)\to\mathcal{CE}(F)$ a.s. as $n\to\infty$; see Proposition 2
in \cite{DicrescenzoLongobardiPROCEEDINGS}. However we can also refer to Proposition
\ref{prop:DL-prop2-generalization} with $w=w_{(1)}$, where
$$w_{(1)}(x):=-x\log x;$$
in fact $C_n(w_{(1)})$ in Proposition \ref{prop:DL-prop2-generalization} coincides with
$\mathcal{CE}(\hat{F}_n)$. It is easy to check that both Conditions \ref{condition-DL-approach}
and \ref{condition} hold for the function $w_{(1)}$.

We can also consider the \emph{fractional generalized cumulative entropy} defined by
$$\mathcal{CE}_\alpha(F):=\frac{1}{\Gamma(\alpha+1)}\int_0^\infty F(z)(-\log F(z))^\alpha dz\ (\mbox{for all}\ \alpha>0)$$
(see eq. (7) in \cite{DicrescenzoKayalMeoli2021}; the generalized cumulative entropy
with $\alpha$ integer was previously defined in \cite{Kayal2016}); note that we recover $\mathcal{CE}(F)$
defined above for $\alpha=1$. In this case we have to consider the function
$$w_{(\alpha)}(x):=\frac{1}{\Gamma(\alpha+1)}x(-\log x)^\alpha;$$
then the function $w=w_{(\alpha)}$ satisfies both Conditions \ref{condition-DL-approach} and \ref{condition}
for $\alpha\geq 1$ while, if $\alpha\in(0,1)$, Condition \ref{condition-DL-approach} holds (with 
$\beta\leq\alpha$) and Condition \ref{condition} fails.

For completeness we also discuss the case of the \emph{fractional cumulative residual 
entropy} (see eq. (5) in \cite{XiongShangZhangCNSNS2019})
$$\mathcal{E}_q(F):=\int_0^\infty (1-F(z))(-\log(1-F(z)))^qdz\ (\mbox{for all}\ q\in [0,1]);$$
note that we recover the cumulative residual entropy defined in \cite{RaoChenVemuriWang} for $q=1$.
In this case we have to consider the function
$$w_{[q]}(x):=(1-x)(-\log(1-x))^q;$$
then, even if we also consider $q>1$, the function $w=w_{[q]}$ satisfies Condition \ref{condition-DL-approach}
with $\beta\in(0,1)$, and does not satisfy Condition \ref{condition} except for $q=0$.

\section*{Funding}
CC and ML are supported by Indam-GNAMPA and by MIUR-PRIN 2017 Project "Stochastic Models
for Complex Systems" (No. 2017JFFHSH). CM and BP are supported by MIUR Excellence Department
Project awarded to the Department of Mathematics, University of Rome Tor Vergata 
(CUP E83C18000100006), by University of Rome Tor Vergata (research program "Beyond Borders", 
project "Asymptotic Methods in Probability" (CUP E89C20000680005)) and by Indam-GNAMPA 
(research project "Stime asintotiche: principi di invarianza e grandi deviazioni").

\section*{Acknowledgements}
We thank Prof. Gao for some discussion on Theorem 4.8 in \cite{GaoZhaoAS2011}.

\section*{Declaration}
The authors declare that they have no conflict of interest.

\end{document}